\documentclass[12pt]{article}
\usepackage{latexsym,amsfonts,amsmath,amsthm,amssymb}

\newcommand{\la}{\lambda}

\newcommand{\Om}{\Omega}
\newcommand{\om}{\omega}
\newcommand{\F}{\cal {F}}
\newcommand{\A}{\cal {A}}
\newcommand{{\bP}}{\bf {P}}

\begin{document}

\date{}

{\bf A formula of total probability with interference term
and the Hilbert space representation of the contextual Kolmogorovian model}

\centerline{Andrei Khrennikov\footnote{Supported in part by the EU Human
Potential Programme, contact HPRN--CT--2002--00279 (Network on
Quantum Probability and Applications) and Profile Math. Modelling
in Physics and Cogn. Sc. of V\"axj\"o University.}}

\centerline{International Center for Mathematical
Modeling}
\centerline{in Physics and Cognitive Sciences,}
\centerline{University of V\"axj\"o, S-35195, Sweden}

\begin{abstract}
We compare the classical Kolmogorov  and quantum
probability models. We show that the gap between these model is not so huge as it was commonly believed.
The main structures of quantum theory (interference of probabilities, Born's rule,
complex probabilistic amplitudes,  Hilbert state space, representation of observables by operators) are present
in a latent form in the  Kolmogorov model. In particular, we obtain ``interference of probabilities'' without to appeal to the Hilbert
space formalism. We interpret ``interference of probabilities'' as a perturbation (by
a $\cos$-term) of the conventional formula of total probability.
Our classical
derivation of quantum probabilistic formalism can stimulate applications of  quantum
methods outside of  microworld : in psychology, biology, economy,...
\end{abstract}



\section{Introduction}

There is a rather common opinion that the quantum model of probability theory (i.e., the calculus on probabilities
based on the complex Hilbert space) differs essentially from the classical (measure-theoretic) Kolmogorov
model  [1], [2]; see, e.g., [3]- [5] for details and discussions. Among distinguishing
features of quantum probability there are typically mentioned:

a) The use of complex amplitudes of probabilities, $\psi(x)$, (wave functions);

b) Born's rule for  probabilities. Probability of the event $B_x$ -- to find a particle at the point $x$ -- is given by
\begin{equation}
\label{BORN}
P_\psi(B_x)= \vert \psi(x) \vert^2.
\end{equation}

c) Interference of probabilities. We present this phenomenon by coupling it to the
formula of total probability. We consider the simplest partition of of the sample space
${\cal A}= \{ A_1, A_2\}.$ Here we have, see, for example, [2]:
\begin{equation}
\label{HYI}
{\bf P}(B\vert C) =\sum {\bf P}(A_j \vert C) {\bf P}(B\vert A_j C).
\end{equation}

However, in the quantum probabilistic formalism there was derived a different formula:
\begin{equation}
\label{INNP}
{\bP}(B\vert C)=\sum {\bP}(A_j\vert C){\bP}(B\vert A_j)+
\end{equation}
\[2 \cos \theta(B\vert {\A},C)\sqrt{{\bP}(A_1\vert C)
{\bP}(B\vert A_1){\bP}(A_2\vert C) {\bP}(B\vert A_2)}\]
where $\theta(B\vert {\A},C)$ is an angle (``phase'') depending on the event $B,$  partition
${\A}$ and the condition $C$ under which the event $B$ occurs. The presence of a new trigonometric term
is interpreted as {\it interference of probabilities}, see, e.g., [6]. In [6] it was emphasized that the presence
of interference of probabilities in quantum formalism is an exhibition of
violation of fundamental laws of classical probability.

d). Representation of physical observables by noncommutative operators in the complex Hilbert space.
(We recall that in the Kolmogorov model there are used random variables -- measurable functions
on the sample space).

The aim of this paper is to show that in fact the gap between quantum model (Dirac-von Neumann [7], [8])
and classical model (Kolmogorov [1]) is not as large as it is commonly believed.\footnote{We do not claim that
all problems are solved. In this paper we do not consider composite systems. Therefore we do not even discuss
such things as  Bell's inequality and quantum nonlocality, see [3], [4] for details.} All mentioned
distinguished features of quantum probability, a)--d), are present in a latent form in the classical
Kolmogorov model.

The crucial point is that all probabilities should be considered as {\it contextual probabilities.}
Here a context $C$ is any complex of conditions, physical, biological, economic, financial. Therefore
it is meaningless to speak about an abstract probability ${\bf P}$ which has no relation to a concrete context.
Any probability should be related to some fixed context $C.$\footnote{
Of course, there is nothing new for probabilists. For example, A. N. Kolmogorov pointed out to the role of
complexes of experimental conditions in defining probability in his famous book [1]
and especially in [9]. Similar views are presented in the books of Gnedenko [10] and Renye [11]. We can also say that
von Mises' frequency probability [12] is contextual: a collective is defined by a complex of experimental
conditions.}

Our main contribution is the contextual probabilistic analysis of the formula of total
probability (\ref{HYI}) and derivation of the ``quantum formula of total probability''
(\ref{INNP}) (which is typically referred to as ``interference of probabilities''). Starting with this formula
(derived in the classical measure-theoretic framework with the Kolmogorov probability space: ${\cal P}= (\Om, \F, {\bP})$) we reproduce other
distinguished features of the quantum probabilistic formalism.

The starting point of our analysis is the contextual interpretation of conditional
probabilities. Typically conditional probability ${\bf P}(A\vert C)$ is interpreted as the probability of occurrence
of the {\it event} $A$ under the condition that the {\it event} $C$ occurred. This interpretation can be called the
{\it event conditioning.} But we would not like to consider conditioning by occurrence of an event.
In general it is impossible to identify, e.g., a collection of equipment
in a laboratory with an event. We consider conditioning by a complex of, e.g., physical conditions $C.$
So our conditioning is {\it conditioning by context and not event.}

An important consequence
of this new interpretation of conditional probabilities ${\bf P}(A\vert C)$ in the Kolmogorov model is
that we are not able to apply Boolean algebra to sets $C$ representing contexts -- complexes of e.g. physical conditions.
For two events, say $C_1$ and $C_2,$  it is always possible to consider the event corresponding to
their simultaneous occurrence. By the Boolean algebra it is realized as $C=C_1 C_2.$ This is a very natural operation
on the algebra of events. But for two contexts it is not always possible to define their simultaneous
realization. Therefore if such contexts are represented by sets $C_1$ and $C_2$  belonging the $\sigma$-algebra
${\cal F}$ of the Kolmogorov space, then by considering the set $C=C_1 C_2$ we cannot be sure that it would represent
a physically meaningful context.

Thus we cannot consider the whole $\sigma$-algebra
${\cal F}$ of the Kolmogorov space as a set-representation of contexts. Depending on a problem under consideration
conditional probabilities ${\bf P}(A\vert C)$ can be considered only for contexts $C$ belonging some special collection ${\cal C} \subset
{\cal F}.$ (An event $A$ is still represented by an arbitrary element of the ${\cal F}).$

We shall show that such a ``cutoff'' of the Kolmogorov $\sigma$-algebra
${\cal F}$ can induce  quantum probabilistic formalism. In such an approach
quantum formalism arises as a special representation of the {\it contextual Kolmogorov model}:
$
{\cal P}_{\rm{cont}}= (\Om, \F\vert {\cal C}, {\bP}))
$
for a special choice of the collection of contexts ${\cal C}.$ \footnote{Finally, we remark that our construction --
the contextual Kolmogorov model -- is very close to Renye's model [11]. Renye also
introduced a special collection of sets, say ${\cal C}_{\rm{REN}},$ representing conditions. But
collections of contexts ${\cal C}$ of our contextual Kolmogorov model do not satisfy conditions of
Renye's model. This gives us the possibility to reproduce  quantum probabilistic formalism that was impossible
to do in Renye's model. The latter model is more general from the measure-theoretic viewpoint. In principle,
we could explore this generality even in our contextual approach. But we shall not do this in the present paper.
We want to show that even the Kolmogorov model contains (in a latent form) main quantum probabilistic structures.
We emphasize again that typically the presence of such structures was considered as an exhibition of non-Kolmogorovness.}

Applying the contextual approach to the formula of total probability (\ref{HYI}), we see that
using of probabilities of the type
${\bf P}(B\vert A_j C),$ i.e., conditioning by ``intersection of contexts'',
in general is meaningless. And we see that
in the ``quantum formula of total probability'' (\ref{INNP}) such probabilities were really excluded from consideration.
Probabilities ${\bf P}(B\vert A C)$
are not defined in the physical framework.
Therefore  in (\ref{INNP}), instead of ${\bf P}(B\vert A C),$
there were considered  ``experimental conditional probabilities''
${\bf P}(B\vert A_j).$
But in general we have the inequality:
\begin{equation}
\label{INNZ}
{\bP}(B\vert C)\not= \sum {\bP}(A_j\vert C){\bP}(B\vert A_j)
\end{equation}
that can be also interpreted as the equality:
\begin{equation}
\label{INNZ1}
{\bP}(B\vert C) = \sum {\bP}(A_j\vert C){\bP}(B\vert A_j)+ \delta(B\vert {\A},C),
\end{equation}
where a perturbation term $\delta(B\vert {\A},C)$ is defined as the difference of the left-hand  and right-hand sides of
(\ref{INNZ}).
In this way, for a special system of contexts ${\cal C}^{\rm{tr}},$ see section 2,
 we obtain the ``quantum formula of total probability'' (\ref{INNP});
and with the aid of this formula we construct a representation of the collection of
contexts ${\cal C}^{\rm{tr}}$ in the unit sphere of the
complex Hilbert space. This is the crucial step to reproduce a)--d) in the classical,
but contextual probabilistic framework.

What are main purposes of such a construction? On one hand, we are able to demystify quantum probability and connect it
in a rather simple way with the classical Kolmogorov model. On the other hand, by reproducing quantum probabilistic
calculus, in particular, ``interference of probabilities'', in the measure-theoretic framework we see that there are no
reasons to restrict applications of this calculus to description of processes in the microworld. By using contextual
approach we can construct the quantum representation for statistical models in any domain of science, for example,
biology, psychology, economics.\footnote{Why can such a representation be fruitful? In our approach the quantum representation
is a {\it projection of the classical probability model.} This is an essential simplification of the classical
probabilistic description. Such a simplified description can be useful for models in that the detailed classical
probabilistic description is extremely complicated, for example,  for applications to cognitive sciences and
psychology and cognitive sciences see [13]--[15], in game theory [16], in financial mathematics [17], in classical theory of disordered
systems [18].}  We remark that the first derivation of the ``quantum formula of total probability'' (\ref{INNP})
 without to appeal to the Hilbert space was done in papers [19], [20] in the {\it von Mises frequency framework}; see also [5]
 for using of the law of large numbers for this purpose.\footnote{Recent years there were also a few attempts to
 use non-Kolmogorovian, but measure-theoretic models to reproduce some predictions of quantum mechanics, see, e.g.,
 [4] and [22].}

\section{Interference formula of total probability}

We consider the conventional formula of total probability (\ref{HYI}) in a special  case.
Let $a$ and $b$ be dichotomous random variables,
$a=a_1, a_2$ and $b=b_1, b_2.$ We have
$$
{\bP}(b=b_i\vert C)=\sum_n {\bP}(a=a_n\vert C) {{\bP}(b=b_i\vert a=a_n, C)} \;.
$$
If a measurement of the variable $a$ disturbs essentially
the context $C,$ then we would not be able to create the context corresponding to nondisturbing measurement
of $a$ under the complex of experimental conditions $C.$ Therefore
we should modify this  formula  and exclude probabilities ${{\bP}(b=b_i\vert a=a_n, C)}.$

The following notion is well known in  measurement theory of quantum mechanics, see [3], [4].
Let us denote by $A_j$ the selection-context with respect to the value $a=a_j$ of the random variable $a$
(for example, in quantum mechanics there are considered momentum-selections: there
are selected all particles with a fixed value of momentum). These contexts $(j=1,2$ in our case)
are represented in the measure-theoretic
approach by sets $A_j=\{\om \in \Om: a(\om)=a_j\}.$ We also introduce the selection-contexts for the $b$-variable. They are
represented by sets
$B_i=\{\om \in \Om: b(\om)=b_i\}.$ We consider partitions ${\cal A}=\{A_1,A_2\}$ and ${\cal B}=\{B_1,B_2\}$
of the sample space $\Omega.$

A set $C$ belonging to $\F$ is said to
be a {\it nondegenerate context} with respect to  the partition ${\A}$ if
${\bP}(A_nC)\not =0$ for all $n.$ We denote the set of all ${\A}-$nondegenerate contexts
 by the symbol  ${\cal C}_{\A, \rm{nd}}.$ The partitions ${\A}$ and ${\cal B}$ are said to be {\it incompatible} if
${\bP}(B_n A_k) \not = 0$ for all $n$ and $k.$ Thus ${\cal B}$ and ${\cal A}$ are incompatible iff every $B_n$ is a nondegenerate context with respect to
${\cal A}$ and vice versa. Random variables $a$ and
$b$ inducing incompatible partitions ${\cal A}$ and ${\cal
B}$  are said to be incompatible. (We remark that we defined incompatibility in purely measure-theoretic framework.)

Everywhere below $a$ and $b$ are incompatible random variables.
Let $B\in {\cal C}_{\A, \rm{nd}}.$ We define a {\it coefficient of interference} of random  variables
$a$  and $b$ by:
\begin{equation}
\label{INTW}
\lambda(B\vert {\cal A}, C)=
\frac{\delta(B\vert {\cal A},C)}{2\sqrt{{\bP}(A_1\vert C) {\bP}(B\vert A_1){\bP}(A_2\vert C) {\bP}(B\vert A_2)}}
\end{equation}
where $\delta(B\vert {\cal A}, C) ={\bf P}(B\vert C) - \sum_{j=1}^{2} {\bf P}(B\vert A_j){\bf P}(A_j\vert C).$
We shall see that the ``perturbed formula of total probability'' (\ref{INNZ1}) has interesting consequences if the
perturbation $\delta$ be represented in the form:
\begin{equation}
\label{INTW1}
\delta(B\vert {\A},C)=\lambda(B\vert {\A},C)\sqrt{{\bP}(A_1\vert C)
{\bP}(B\vert A_1){\bP}(A_2\vert C) {\bP}(B\vert A_2)}
\end{equation}
We set
$$
{\cal C}^{\rm{tr}} =\{ C\in {\cal C}_{\A, \rm{nd}}: \vert\lambda(B_j\vert {\cal A}, C)\vert\leq 1\}
$$
We call elements of ${\cal C}^{\rm{tr}}$ trigonometric contexts. We consider the contextual Kolmogorov model
with this collection of contexts:
\begin{equation}
\label{INTW1A}
{\cal P}_{\rm{cont,tr}}= (\Om, \F\vert {\cal C}^{\rm{tr}}, {\bP})
\end{equation}
We remark that in general the system of sets ${\cal C}^{\rm{tr}}$
is not an algebra: $C_1, C_2 \in {\cal C}^{\rm{tr}}$ does not imply
that $C=C_1 C_2 \in {\cal C}^{\rm{tr}}.$
Our main result can be formulated in the form of the following theorem (which will be proved in a few steps):

\medskip

{\bf Theorem 2.1.} {\it The ``quantum formula of total probability''
(\ref{INNP}) can be derived in the Kolmogorov probability framework. On the basis of this formula
we can construct a  map  from the set of trigonometric contexts ${\cal C}^{\rm{tr}}$
into the unit sphere $S$ of the complex Hilbert space $H$ (space of complex amplitudes).
Such a map is determined by a pair $a, b$ of incompatible random variables (reference variables) that are represented
by noncommutative operators $\hat{a}, \hat{b}$.
Unitarity of the matrix $V^{b\vert a}$ of transition from the basis $\{e^a_i\}$
to the basic $\{e_i^b\}$ (these bases correspond to random variables $a$ and
$b$) is equivalent
to Born's rule for both reference variables. This construction can be realized only for a double stochastic matrix
of transition probabilities.}

\medskip

First by using the relation (\ref{INTW1}) we see that
the ``perturbed formula of total probability'' (\ref{INNZ1})  can be written as:
\begin{equation}
\label{INN}
{\bP}(B\vert C)=\sum {\bP}(A_j\vert C){\bP}(B\vert A_j)+
\end{equation}
\[2 \lambda(B\vert {\A},C)\sqrt{{\bP}(A_1\vert C)
{\bP}(B\vert A_1){\bP}(A_2\vert C) {\bP}(B\vert A_2)}\]

1). Suppose that the interference coefficients
$\vert \la(B\vert {\A},C)\vert\leq 1$ for every $B\in {\cal B}.$
We introduce
new statistical parameters $\theta(B\vert {\A},C)\in
[0,2 \pi]$ and represent the coefficients in the
trigonometric form:
$
\la(B\vert {\A},C)=\cos \theta (B\vert {\A},C).
$
Parameters $\theta(B\vert {\A},C)$ are said to be {\it{relative phases}} of
an event $B$ with respect to the partition ${\A}$
in the context $C$.
In this case the ``perturbed formula of total probability'' given in the form (\ref{INN}) coincides with
the ``quantum formula of total probability'' (\ref{INNP}).
\footnote{This is nothing other than the famous formula of interference of
probabilities. Typically this formula is derived by using the Hilbert
space (unitary) transformation corresponding to the transition from
one orthonormal basis to another and Born's probability postulate.
The orthonormal basis under quantum consideration consist of eigenvectors of
operators (noncommutative) corresponding to quantum physical observables
$a$ and $b.$}

2). Suppose that
$\vert \la(B\vert {\A},C)\vert\geq  1$ for every $B\in {\cal B}.$
We set $\theta(B\vert {\A},C)\in
(-\infty ,+ \infty)$ and represent the coefficients  in the
hyperbolic form:
$
\la(B\vert {\A},C)=\pm \cosh \theta(B\vert {\A},C).
$
In this case (\ref{INN}) has the form of ``hyperbolic interference of probabilities''\begin{equation}
\label{TNC1}
{\bP}(B\vert C)=\sum {\bP}(A_j\vert C){\bP}(B\vert A_j)\pm
\end{equation}
\[2 \cosh \theta(B\vert {\A},C)\sqrt{{\bP}(A_1\vert C)
{\bP}(B\vert A_1){\bP}(A_2\vert C) {\bP}(B\vert A_2)}\]

In this paper we shall concentrate our considerations on the first case.\footnote{We just mention that
in the second case we can obtain a representation of the contextual Kolmogorov model
$
{\cal P}_{\rm{cont, hyp}}= (\Om, \F\vert {\cal C}^{\rm{hyp}}, {\bP}),
$
where
$
{\cal C}^{\rm{hyp}} =\{ C\in {\cal C}_{\A, \rm{nd}}: \vert\lambda(B_j\vert {\cal A}, C)\vert \geq 1\},
$
in so called hyperbolic Hilbert space: a Hilbert module over the two dimensional Clifford algebra
(i.e., the commutative algebra with basis $e_1=1$ and $e_2=j$, where $j^2=+1,$ see [21] for details).
Therefore  it is impossible to represent the whole Kolmogorov $\sigma$-algebra  ${\cal F}$in the complex Hilbert
space. Moreover, ${\cal C}^{\rm{tr}} \cup {\cal C}^{\rm{hyp}}$
is a proper sub-system of ${\cal F}.$ For example, there exist mixed hyper-trigonometric contexts: one
$\lambda \leq 1$ and another $\lambda \geq 1$.  There also exist degenerate
contexts $C$ for that interference coefficients are not defined at all.}

Everywhere below $B=B_x, x=b_1,b_2,$ and we shall often use the symbols  $ \lambda(b=x\vert a, C)$
instead of $\lambda(B_x\vert {\cal A}, C).$

\section{Extraction of complex probability amplitudes and Born's rule from the Kolmogorov model}

We recall that we study  the case of incompatible dichotomous random
variables $a=a_1, a_2, b=b_1, b_2.$ This pair of variables will be fixed.
We call such variables {\bf reference variables.} For each fixed pair $a, b$ of reference variables
we construct a representation of the contextual Kolmogorov model
${\cal P}_{\rm{cont,tr}}= (\Om, \F\vert {\cal C}^{\rm{tr}}, {\bP}))$
in the complex Hilbert space. We set $Y=\{a_1, a_2\}, X=\{b_1,
b_2\}$ (``spectra'' of random variables $a$ and $b).$
Let $C\in {\cal C}^{\rm{tr}}.$
 We set
 $$
 p_C^a(y)={\bP}(a=y\vert C), p_C^b(x)={\bP}(b=x\vert C), p(x\vert y)={\bP}(b=x\vert a=y),
 $$
$x \in X, y \in Y.$
The  formula  (\ref{INNP}) can be written as
\begin{equation}
\label{Two}
p_c^b(x)=\sum_{y \in Y}p_C^a(y) p(x\vert y) + 2\cos \theta_C(x)\sqrt{\Pi_{y \in
Y}p_C^a(y) p(x\vert y)}\;,
\end{equation}
where $\theta_C(x)=\theta(b= x\vert a, C)= \pm \arccos \lambda(b=x\vert a, C), x \in X.$
Here\\
$
\delta(b=x\vert a, C)=p_c^b(x)-\sum_{y \in Y} p_C^a(y)p(x\vert y)
$
and
$\lambda(b=x\vert a, C)  =\frac{\delta(b=x\vert a, C)}{2\sqrt{\Pi_{y\in Y}p_C^a(y)p(x\vert y)}} .$
By using the elementary formula:
$
D=A+B+2\sqrt{AB}\cos \theta=\vert \sqrt{A}+e^{i \theta}\sqrt{B}|^2,
$
for $A, B > 0, \theta\in [0,2 \pi],$
we can represent the probability $p_C^b(x)$ as the square of the complex amplitude (Born's rule):
\begin{equation}
\label{Born}
p_C^b(x)=\vert\varphi_C(x)\vert^2 \;.
\end{equation}
We set
\begin{equation}
\label{EX1}
\varphi(x)\equiv \varphi_C(x)=\sqrt{p_C^a(a_1)p(x\vert a_1)} + e^{i \theta_C(x)} \sqrt{p_C^a(a_2)p(x\vert a_2)} \;.
\end{equation}

It is important to underline that since  for each $x\in X$ phases $\theta_C(x)$ can be
chosen in two ways (by choosing signs + or -) a representation of contexts by complex
amplitudes is not uniquely defined.\footnote{To fix a representation of the contextual
Kolmogorov space ${\cal P}_{\rm{cont, tr}}$ we should fix phases. We shall see that to obtain a
``good representation'' we should choose phases in a special way.}

We denote the space of functions: $\varphi: X\to {\bf C},$ where ${\bf C}$ is the field
of complex numbers, by the symbol
$E=\Phi(X, {\bf C}).$ Since $X= \{b_1, b_2 \},$ the $E$ is the two dimensional
complex linear space. Dirac's $\delta-$functions $\{ \delta(b_1-x), \delta(b_2-x)\}$
form the canonical basis in this space. We shall see (Proposition 5.1)
that under natural assumption
on the matrices of transition probabilities $\varphi_{B_j}(x)= \delta(b_j - x).$
For each $\varphi \in E$ we have
$\varphi(x)=\varphi(b_1) \delta(b_1-x) + \varphi(b_2) \delta(b_2-x).$
By using the representation (\ref{EX1}) we construct the map
\begin{equation}
\label{MAP}
J^{b\vert a}:{\cal C}^{\rm{tr}} \to \Phi(X, {\bf C})
\end{equation}
The $J^{b\vert a}$ maps contexts (complexes of, e.g., physical conditions) into complex
amplitudes. The representation ({\ref{Born}}) of probability as the square of the
absolute value of the complex $(b\vert a)-$amplitude is nothing other than the
famous {\bf Born rule.} The complex amplitude $\varphi_C(x)$ can be called a {\it wave function} of the
complex of physical conditions (context $C)$  or a  pure {\it state.}
 We set
$e_x^b(\cdot)=\delta(x- \cdot).$
The Born's rule for complex amplitudes (\ref{Born}) can be rewritten in the following form:
\begin{equation}
\label{BH}
p_C^b(x)=\vert(\varphi_C, e_x^b)\vert^2 \;,
\end{equation}
where the scalar product in the space $E=\Phi(X, {\bf C})$ is defined by the
standard formula:
$(\varphi, \psi) = \sum_{x\in X} \varphi(x)\bar \psi(x).$
The system of functions $\{e_x^b\}_{x\in X}$ is an orthonormal basis in the
Hilbert space $H=(E, (\cdot, \cdot))$
Let $X \subset {\bf R},$ where ${\bf R}$ is the field of real numbers. By using the Hilbert space representation ({\ref{BH}}) of
the Born's rule  we obtain  the Hilbert space representation of the
expectation of the (Kolmogorovian) random variable $b$:
\begin{equation}
\label{BI1}
E (b\vert C)= \sum_{x\in X}xp_C^b(x)=\sum_{x\in X}x\vert\varphi_C(x)\vert^2=
\sum_{x\in X}x (\varphi_C, e_x^b) \overline{(\varphi_C, e_x^b)}=
(\hat b \varphi_C, \varphi_C) \;,
\end{equation}
where  the  (self-adjoint) operator $\hat b: E \to E$ is determined by its
eigenvectors: $\hat b e_x^b=x e^b_x, x\in X.$
This is the  multiplication operator in the space of complex functions $\Phi(X,{\bf C}):$
$
\hat{b} \varphi(x) = x \varphi(x)
$
 By (\ref{BI1}) the  conditional expectation of the Kolmogorovian
random variable $b$ is represented
with the aid of the self-adjoint operator $\hat b.$  Therefore it is natural to
represent this random variable (in the Hilbert space model)  by the operator $\hat b.$
We shall use the following notations:
\begin{equation}
\label{KOE}
u_j^a=\sqrt{p_C^a(a_j)}, u_j^b=\sqrt{p_C^b(b_j)}, p_{ij}=p(b_j\vert a_i), u_{ij}=\sqrt{p_{ij}},
\theta_j=\theta_C(b_j).
\end{equation}
We remark that the coefficients $u_j^a, u_j^b$ depend on a context $C;$ so
$u_j^a=u_j^a(C), u_j^b=u_j^b(C).$
We also consider the {\it matrix of transition
probabilities} ${\bf P}^{b\vert a}=(p_{ij}).$
It is always a {\it stochastic matrix.}\footnote{So $p_{i1}+p_{i2}=1, i=1,2.$}
We have
$
\varphi_C=v_1^b e_1^b + v_2^b e_2^b, \;\mbox{where}\;\;
v_j^b=u_1^a u_{1j}  + u_2^a u_{2j} e^{i \theta_j}\;.
$
Hence
\begin{equation}
\label{BI}
p_C^b(b_j) =\vert v_j^b \vert^2 = \vert u_1^a u_{1j}  + u_2^a u_{2j} e^{i \theta_j} \vert^2.
\end{equation}

This is the {\it interference representation of probabilities} that is used,
e.g., in quantum formalism. We recall that we obtained (\ref{BI}) starting
with the interference formula of total probability, (\ref{Two}).

We would like to have Born's rule not only for the $b$-variable,
but also for the $a$-variable. As we shall see, we cannot be lucky in the general case.
Starting from two arbitrary incompatible (Kolmogorovian) random variables
$a$ and $b$ we obtained a complex linear space representation of the
probabilistic model which is essentially more general than the standard quantum
representation. In our (more general) linear representation the ``conjugate
variable'' $a$ need not be represented by a symmetric operator (matrix) in
the Hilbert space $H$ generated by the $b$. We recall that in QM both reference
variables (the position and the momentum) are represented in the same Hilbert
space.

For any context $C_0,$ we can represent the corresponding wave function
$\varphi=\varphi_{C_0}$ in the form:
\begin{equation}
\label{0}
\varphi=u_1^a e_1^a + u_2^a e_2^a,
\end{equation}
where
\begin{equation}
\label{Bas}
e_1^a= (u_{11}, \; \; u_{12}) ,\; \;
e_2^a= (e^{i \theta_1} u_{21}, \; \; e^{i \theta_2} u_{22})
\end{equation}
Here $\{e_i^a\}$ is a system of vectors in $E$ corresponding to the $a$-observable.
We suppose that vectors $\{e_i^a\}$ are lineary independent, so $\{e_i^a\}$
is a basis in $E.$ We have:
$e_1^a=v_{11} e_1^b + v_{12} e_2^b, \; \; \; e_2^a=v_{21} e_1^b + v_{22} e_2^b.$
Here $V=(v_{ij})$
is the matrix corresponding to the transformation of complex amplitudes:
$v_{11}=u_{11}, v_{21}=u_{21}$ and $v_{12}=e^{i \theta_1} u_{21}, v_{22}=
e^{i \theta_2} u_{22}.$
We would like to find a class of matrixes $V$ such that
Born's rule (in the Hilbert space form), see (\ref{BH}), holds true also in the $a-$basis:
\begin{equation}
\label{BBR}
p_C^a(a_j)=\vert(\varphi, e_j^a)\vert^2 \; .
\end{equation}
By (\ref{0}) we have the
Born's rule (\ref{BBR}) iff $\{e_i^a\}$ was an {\it orthonormal
basis,} i.e., the  $V$ was a {\it unitary} matrix.

Since we study the
two-dimensional case (i.e., dichotomous random variables), $V\equiv
V^{b\vert a}$ is unitary iff the matrix of transition probabilities ${\bf
P}^{b\vert a}$ is {\bf double stochastic} and $e^{i\theta_1}=-e^{i\theta_2}$ or
\begin{equation}
\label{MAR}
\theta_{C_0}(b_1) - \theta_{C_0}(b_2)=\pi \mod 2\pi
\end{equation}

We recall that a matrix is double stochastic if it is stochastic, i.e.,
$p_{j1} + p_{j2}=1,$ and, moreover,
$p_{1j} + p_{2j}=1, j=1,2.$
Any matrix of transition probabilities is stochastic,
but in general it is not double stochastic.
We remark that the constraints (\ref{MAR}) on phases and the double stochasticity constraint
are not independent:

{\bf Lemma 3.1.} {\it Let the
matrix of transition probabilities ${\bf P}^{b\vert a}$ be double stochastic.
Then:
\begin{equation}
\label{SW}
\cos \theta_C (b_2)=-\cos \theta_C (b_1)
\end{equation}
for any context $C\in {\cal C}^{\rm{tr}}.$}

By Lemma 3.1 we have two different possibilities to choose phases:
\[\theta_{C_0}(b_1) + \theta_{C_0}(b_2) = \pi  \;\rm{or} \;
\theta_{C_0}(b_1) - \theta_{C_0}(b_2) = \pi \mod 2\pi\]
By (\ref{MAR}) to obtain the Born's rule for the $a$-variable we should choose phases $\theta_{C_0}(b_i), i=1,2,$ in such a way that
\begin{equation}
 \label{MAR0}
 \theta_{C_0}(b_2)=\theta_{C_0}(b_1) + \pi.
\end{equation}
If $\theta_{C_0}(b_1)\in [0, \pi]$ then $\theta_{C_0}(b_2)\in [\pi, 2\pi]$ and vice versa.
Lemma 3.1 is very important since by it (in the case when reference observables are chosen
in such way that the matrix of transition probabilities is double stochastic)
we can always choose $\theta_{C_0}(b_j), j=1,2,$ to satisfy (\ref{MAR0}).

The delicate feature of the presented construction of the $a$-representation is that the basis $e_j^a$ depends on the context $C_0: e_j^a=e_j^a(C_0).$ And the Born's rule, in fact, has the form:
$p_{C_0}^a(a_j)=|(\varphi_{C_0}, e_j^a(C_0))|^2.$
We would like to use (as in the conventional quantum formalism)
one fixed $a$-basis for all contexts $C\in {\cal C}^{\rm{tr}}.$ We may try to use for all contexts
$C\in {\cal C}^{\rm{tr}}$ the basis $e_j^a\equiv e_j^a(C_0)$ corresponding to one fixed context $C_0.$ We shall
see that this is  really the  fruitful  strategy.

\medskip

{\bf Lemma 3.2}
{\it {Let the matrix of transition probabilities ${\bf P}^{b\vert a}$ be double
stochastic and let for any context $C\in {\cal C}^{\rm{tr}}$ phases $\theta_C(b_j)$ be chosen as
\begin{equation}
\label{P}
\theta_C(b_2)=\theta_C(b_1) + \pi \mod 2\pi.
\end{equation}
Then for any context $C\in {\cal C}^{\rm{tr}}$ we have the Born's rule for the basis $e_j^a\equiv e_j^a(C_0)$
constructed for a fixed context}} $C_0:$
\begin{equation}
\label{MAR2}
p_C^a(a_j)=|(\varphi_C, e_j^a)|^2
\end{equation}

{\bf Proof.} Let $C_0$ be some fixed context. We take the basic
$\{ e_j^a(C_0)\}$ (and the matrix $V(C_0))$ corresponding to this context.
For any $C \in {\cal C}^{\rm{tr}},$ we would like to represent the wave function $\phi_C$ as
$\phi_C= v_1^a(C) e_1^a(C_0) + v_2^a(C) e_2^a(C_0),\;\; \mbox{where}\;\;\; \vert v_j^a(C) \vert^2= p_C^a(a_j).
$
It is clear that, for any $C \in {\cal C}^{\rm{tr}},$ we can represent the wave function as
$$
\phi_C(b_1) = u_1^a(C) v_{11}(C_0) + e^{i[\theta_C(b_1)- \theta_{C_0}(b_1)]} u_2^a(C) v_{12}(C_0)
$$
$$
\phi_C(b_2) = u_1^a(C) v_{21}(C_0) + e^{i [\theta_C(b_2)- \theta_{C_0}(b_2)]} u_2^a(C) v_{22}(C_0)
$$
Thus we should have:
$\theta_C(b_1)- \theta_{C_0}(b_1)=
\theta_C(b_2)- \theta_{C_0}(b_2)   \mod  2\pi.$
for any pair of contexts $C_0$ and $C_1.$
By using the relations (\ref{P}) between phases
$\theta_C(b_1), \theta_C(b_2)$ and $\theta_{C_0}(b_1), \theta_{C_0}(b_2)$ we obtain:
$
\theta_C(b_2)-\theta_{C_0}(b_2)=(\theta_C(b_1) + \pi-\theta_{C_0}(b_1)-\pi)=
\theta_C(b_1)-\theta_{C_0}(b_1) \mod 2\pi.
$

\medskip

The constraint (\ref{P}) essentially restricted the class of complex amplitudes which can
be used to represent a context $C\in {\cal C}^{\rm{tr}}$. Any $C$ can be represented only by
two amplitudes $\varphi(x)$ and $\bar{\varphi}(x)$ corresponding to the two possible
choices of $\theta_C(b_1)$ (in $[0, \pi]$ or $(\pi, 2\pi$)).

By Lemma 3.2 we obtain the following part of the Theorem 2.1:
{\it We can construct the Hilbert space representation of
the contextual Kolmogorov model ${\cal P}_{\rm{cont,tr}}$ such that the Born's rule holds true for both reference
variables iff the matrix of transition probabilities ${\bf P}^{b\vert a}$ is double stochastic.}

If ${\bf P}^{b\vert a}$ is double stochastic we have a quantum-like representation not only for
the conditional expectation of the variable $b,$ see (\ref{BI1}), but also for the variable $a:$
\begin{equation}
\label{EXP}
E(a\vert C) = \sum_{y\in Y} y p_C^a(y)= \sum_{y\in Y} y \vert(\varphi_C, e_y^a)\vert^2
=(\hat {a}\varphi_C, \varphi_C) \;,
\end{equation}
where the self-adjoint operator (symmetric matrix)
$\hat{a} :E\to E$ is determined by its eigenvectors: $\hat{a}
e_j^a=a_j e_j^a.$  By (\ref{EXP}) it is natural to represent the random variable
$a$ by the operator $\hat{a}.$

Let us denote the unit sphere in the Hilbert space $E=\Phi(X, {\bf C})$
by the symbol $S.$ The map $J^{b\vert a}:{\cal C}^{\rm{tr}}\to S$ need not be a surjection (injection).
In general the set of (pure) states corresponding to a contextual Kolmogorov
space
$$
S_{{\cal C}^{\rm{tr}}}\equiv S^{b\vert a}_{{\cal C}^{\rm{tr}}}
=J^{b\vert a}({\cal C}^{\rm{tr}})
$$
is just a proper subset of the sphere $S.$ The structure of the set of pure states
$S_{{\cal C}^{\rm{tr}}}$ is determined by the Kolmogorov space and the reference variables
$a$ and $b.$

\section{Noncommutative operator-representation of Kolmogorovian random variables}

Let the matrix of transition probabilities
${\bf P}^{b\vert a}$ be {\bf double stochastic.}
We consider in this section the case of real valued random variables. Here
spectra  of random variables $b$ and $a$ are subsets of ${\bf R}.$
We set $q_1= \sqrt{p_{11}}=\sqrt{p_{22}}$ and $q_2= \sqrt{p_{12}}= \sqrt{p_{21}}.$
Thus the vectors of the $a$-basis, see (\ref{Bas}), have the following form:
$
e_1^a= (q_1, q_2), \; \; e_2^a= (e^{i \theta_1} q_2, e^{i \theta_2} q_1)\;.
$
Since $\theta_2 = \theta_1+ \pi,$ we get $e_2^a= e^{i \theta_2} (- q_2,  q_1).$
We now find matrices of operators $\hat{a}$ and $\hat{b}$ in the $b$-representation. The latter one
is diagonal. For $\hat{a}$ we have:
$\hat{a}= V \rm{diag}(a_1, a_2) V^\star,$ where $v_{11}=v_{22}=q_1, v_{21}=-v_{12}=q_2.$ Thus
$
a_{11}= a_1q_1^2 +a_2 q_2^2, \; a_{22}= a_1q_2^2 +a_2 q_1^2,\;
a_{12} = a_{21}= (a_1-a_2) q_1 q_2 .
$
Hence
$
[\hat{b}, \hat{a}] = \hat{m},
$
where $m_{11}=m_{22}=0$ and $m_{12}=- m_{21}= (a_1- a_2) (b_2-b_1) q_1 q_2.$
Since $a_1\not= a_2, b_1\not= b_2$ and $q_j\not=0,$ we have $\hat{m}\not=0.$

\section{The role of simultaneous double stochasticity of ${\bf P}^{b\vert a}$ and
${\bf P}^{a\vert b}$}

Starting with the $b$-representation -- complex amplitudes $\phi_C(x)$ defined on
the spectrum (range of values) of a random variable $b$ -- we
constructed the $a$-representation. This construction is natural (i.e., it  produces the
Born's probability rule) only when the  ${\bf P}^{b\vert a}$ is double stochastic.
We would like to have a symmetric model. So by starting with
the $a$-representation -- complex amplitudes $\phi_C(y)$ defined on
the spectrum (range of values) of a random variable $a$ -- we would like
to construct the natural $b$-representation. Thus both matrices of transition
probabilities ${\bf P}^{b\vert a}$ and ${\bf P}^{a\vert b}$ should be double stochastic.

{\bf Theorem 5.1.} {\it Let the matrix ${\bf P}^{b\vert a}$ be double stochastic. The contexts
$B_1, B_2$ belong to ${\cal C}^{\rm{tr}}$ iff the matrix ${\bf P}^{a\vert b}$ is double stochastic.}

{\bf Lemma 5.1.} {\it Both matrices of transition probabilities ${\bf P}^{b\vert a}$ and ${\bf P}^{a\vert b}$
are double stochastic iff the transition probabilities are symmetric, i.e.,
\begin{equation}
\label{SYM}
p(b_i\vert a_j)=p(a_j\vert b_i), i, j=1,2 .
\end{equation}
This is equivalent that random variables $a$ and $b$ have the uniform probability distribution:
$p^a(a_i)=p^b(b_i)=1/2, i=1,2.$}

This Lemma has important physical consequences. A natural (Bornian) Hilbert space representation
of contexts can be constructed only on the basis of a pair of (incompatible) uniformly distributed
random variables.

{\bf Lemma 5.2.}
{\it{Let both matrices ${\bf P}^{b\vert a}$ and ${\bf P}^{a\vert b}$ be double stochastic. Then}}
\begin{equation}
\label{L}
\la (B_i\vert a, B_i) =1 .
\end{equation}

{\bf Proposition 5.1.}
{\it{Let both matrices of transition probabilities ${\bf P}^{b\vert a}$ and ${\bf P}^{a\vert b}$ be double stochastic. Then}}
$$
J^{b\vert a}(B_j)(x)=\delta(b_j-x), x  \in X, \;\; \mbox{and} \; \; \; J^{a\vert b}(A_j)(y)=\delta(a_j-y), y \in Y.
$$
\noindent

Thus in the case when both matrices of transition probabilities
${\bf P}^{a\vert b}$ and ${\bf P}^{b\vert a}$ are double stochastic (i.e., both reference
variables $a$ and $b$ are uniformly distributed) the Born's rule has the form:
$p_C^b(x) = \vert (\phi_C, \phi_{B_x})\vert^2.$

\section{Complex amplitudes of probabilities in the case of multivalued reference variables}

The general case of random variables taking $n\geq 2$ different values can be (inductively) reduced
to the case of dichotomous random variables. We consider two incompatible
random variables taking $n$ values: $b=b_1, \ldots, b_n$ and
$a=a_1, \ldots, a_n.$
We start with some evident generalizations of results presented in section 2.

{\bf Lemma 6.1.} {\it{Let $B, C, D_1, D_2 \in {\cal F}, {\bf P}(C)\ne 0$ and $D_1\cap D_2=\emptyset.$ Then}}
\begin{equation}
\label{F1}
{\bf P}(B(D_1 \cup D_2)\vert C)={\bf P}(BD_1\vert C)+{\bf P}(BD_2\vert C)
\end{equation}

{\bf Proposition 6.1.} (The formula of total probability) {\it{Let conditions of Lemma 6.1 hold
true and let ${\bf P}(D_jC)\ne 0.$ Then}}
\begin{equation}
\label{F2}
{\bf P}(B(D_1 \cup D_2)\vert C)={\bf P}(B\vert D_1C){\bf P}(D_1\vert C)+{\bf P}(B\vert D_2C){\bf P}(D_2\vert C)
\end{equation}

{\bf Proposition 6.2.} (Contextual formula of total probability){\it{Let conditions of Proposition 6.1 hold true
and let ${\bf P}(BD_j)\ne 0, j=1,2.$ Then
\begin{equation}
\label{F3}
{\bf P}(B(D_1 \cup D_2)\vert C)={\bf P}(B\vert D_1){\bf P}(D_1\vert C)+{\bf P}(B\vert D_2){\bf P}(D_2\vert C)+
\end{equation}
\[2\lambda(B\vert \{D_1, D_2\}, C) \sqrt{{\bf P}(B\vert D_1){\bf P}(D_1\vert C){\bf P}(B\vert D_2){\bf P}(D_2\vert C)},\]
where the ``interference coefficient''
\begin{equation}
\label{LA}
\lambda(B\vert \{D_1, D_2\}, C)=\frac{\delta(B\vert \{D_1, D_2\}, C)}{2\sqrt{{\bf P}(B\vert D_1){\bf P}(D_1\vert C)
{\bf P}(B\vert D_2){\bf P}(D_2\vert C)}}
\end{equation}
and }
$\delta(B\vert \{D_1, D_2\}, C)={\bf P}(B(D_1 \cup D_2)\vert C) - \sum_{j=1}^{2} {\bf P}(B\vert D_j){\bf P}(D_j\vert C)$
\[=\sum_{j=1}^{2} {\bf P}(D_j\vert C) ({\bf P}(B\vert D_j C)-{\bf P}(B\vert D_j))\]}

\medskip

We remark that if ${\cal D}=\{ D_1, D_2\}$ is a partition of the sample space, then the formula (\ref{F3})
coincides with the interference formula of total
probability, see section 2.

In the construction of a Hilbert space representation of contexts for multivalued observables there
will be used the following combination of formulas (\ref{F1}) and (\ref{F3}).

{\bf Lemma 6.2.} {\it Let conditions of Lemma 6.1 hold true and let ${\bf P}(BD_1),$ \\${\bf P}(CD_1)$
and ${\bf P}(BD_2 C)$ be strictly positive. Then
\begin{equation}
\label{F5}
{\bf P}(B(D_1 \cup D_2)\vert C)={\bf P}(B\vert D_1){\bf P}(D_1\vert C)+ {\bf P}(BD_2\vert C)
\end{equation}
\[+ 2 \mu(B\vert \{D_1, D_2\}, C) \sqrt{{\bf P}(B\vert D_1){\bf P}(D_1\vert C){\bf P}(BD_2\vert C)}\]
where} $\mu(B\vert \{D_1, D_2\}, C)=\frac{{\bf P}(B(D_1 \cup D_2)\vert C)-{\bf P}(B\vert D_1){\bf P}(D_1\vert C)-{\bf P}(BD_2\vert C)}{2\sqrt{{\bf P}(B\vert D_1){\bf P}(D_1\vert C){\bf P}(BD_2\vert C)}
}$

\medskip

Suppose that coefficients  $\mu$ and $\lambda$ are bounded
by 1. Then we can represent them in the trigonometric form:
\[
\lambda(B\vert \{D_1, D_2\}, C)=\cos \theta (B\vert \{D_1, D_2\}, C)
\]
\[
\mu(B\vert \{D_1, D_2\}, C) =\cos \gamma (B\vert \{D_1, D_2\}, C)
\]

By inserting these cos-expressions in (\ref{F3}) and (\ref{F5}) we obtain trigonometric transformations of
probabilities. We have (by Lemma 6.2):
$$
{\bf P}(B_x\vert C)={\bf P}(B_x(A_1 \cup \ldots \cup A_n)\vert C)
$$
$$
={\bf P}(B_x\vert A_1){\bf P}(A_1\vert C)
+ {\bf P}(B_x(A_2 \cup \ldots \cup A_n)\vert C)
$$
$$
+ 2\mu(B_x\vert \{A_1, A_2 \cup \ldots \cup A_n\}, C) \sqrt{ {\bf P}(B_x\vert A_1){\bf P}(A_1\vert C){\bf P}(B_x(A_2 \cup \ldots \cup A_n)\vert C)},
$$
where
$
\mu(B_x\vert \{A_1, A_2 \cup \ldots \cup A_n\}, C)
$
$$
=\frac{{\bf P}(B_x(A_1 \cup \ldots \cup A_n)\vert C) - {\bf P}(B_x\vert A_1){\bf P}(A_1\vert C)-{\bf P}(B_x(A_2 \cup \ldots \cup A_n)\vert C)}
{2\sqrt{{\bf P}(B_x\vert A_1){\bf P}(A_1\vert C){\bf P}(B_x(A_2 \cup \ldots \cup A_n)\vert C))}} .
$$
Suppose that the coefficients  are relatively small for all $x \in X:$
$|\mu(B_x\vert \{A_1, A_2 \cup \ldots \cup A_n\}, C)|\leq 1.$
Then we can represent these coefficients as
$
\mu(B_x\vert \{A_1, A_2 \cup \ldots \cup A_n\}, C) = \cos \gamma (B_x\vert \{A_1, A_2 \cup \ldots \cup A_n\}, C).
$
Thus the probability ${\bf P}(B_x\vert C)\equiv {\bf P}(B_x(A_1 \cup \ldots \cup A_n)\vert C)$
can be represented as the square of the absolute value of the complex amplitude:
$$
\varphi_C(x)\equiv \varphi_C^{(1)}(x)=\sqrt{{\bf P}(B_x\vert A_1){\bf P}(A_1\vert C)} +
e^{i\gamma_C^{(1)}(x)} \sqrt{{\bf P}(B_x(A_2 \cup \ldots \cup A_n)\vert C)},
$$
where the phase $\gamma_C^{(1)} (x) \equiv \gamma(B_x\vert \{A_1, A_2 \cup \ldots \cup A_n\}, C).$
In the same way the probability in the second summand can be represented as:
$$
{\bf P}(B_x(A_2 \cup \ldots \cup A_n)\vert C)=
{\bf P}(B_x\vert A_2){\bf P}(A_2\vert C)+ {\bf P}(B_x(A_3 \cup \ldots \cup A_n)\vert C) +
$$
$$
2\mu(B_x\vert \{A_2, A_3 \cup \ldots \cup A_n\}, C)
\sqrt{{\bf P}(B_x\vert A_2){\bf P}(A_2\vert C) {\bf P}(B_x(A_3 \cup \ldots \cup A_n)\vert C)},
$$
where
$$
\mu(B_x\vert \{A_2, A_3 \cup \ldots \cup A_n\}, C)
$$
$$
=\frac{{\bf P}(B_x(A_2 \cup \ldots \cup A_n)\vert C)-{\bf P}(B_x\vert A_2){\bf P}(A_2\vert C)-{\bf P}(B_x(A_3 \cup \ldots \cup A_n)\vert C)}{2\sqrt{{\bf P}(B_x\vert A_2){\bf P}(A_2\vert C){\bf P}(B_x(A_3 \cup \ldots \cup A_n)\vert C)}}.
$$

By supposing that these coefficients of statistical disturbance are bounded by 1 we represent the probability as the square of the absolute value of the complex amplitude:
$$
\varphi_C^{(2)}(x)=\sqrt{{\bf P}(B_x\vert A_2){\bf P}(A_2\vert C)} + e^{i \gamma_C^{(2)}(x)} \sqrt{{\bf P}(B_x(A_3 \cup \ldots \cup A_n)\vert C)},
$$
where $\gamma_C^{(2)}(x)=\pm \arccos \mu (B_x\vert \{A_2, A_3, \cup \ldots \cup A_n\}, C).$ On the $jth$ step we represent ${\bf P}(B_x(A_j \cup \ldots \cup A_n)\vert C)$ as the square of the absolute value of the complex amplitude
\[\varphi_C^{(j)}(x)=\sqrt{{\bf P}(B_x\vert A_j){\bf P}(A_j\vert C)} + e^{i \gamma_C^{(j)}(x)} {\sqrt{{\bf P}(B_x(A_{j+1}\cup \ldots \cup A_n)\vert C)}},\]
where $\gamma_C^{(j)}(x)$ is the phase of the coefficient
$$
\mu (B_x\vert \{A_j, A_{j+1} \cup \ldots \cup A_n\}, C)
$$
$$
=\frac{{\bf P}(B_x(A_j \cup \ldots \cup A_n)\vert C)-{\bf P}(B_x\vert A_j) {\bf P}(A_j\vert C)-{\bf P}(B_x (A_{j+1}\cup \ldots \cup A_n)\vert C)}{2\sqrt{{\bf P}(B_x\vert A_j) {\bf P}(A_j\vert C) {\bf P}(B_x (A_{j+1} \cup \ldots \cup A_n)\vert C)}}.$$

It is supposed that at each step we obtain coefficients $|\mu|$ bounded by 1.
At the step $j=n-1$ we should represent the probability ${\bf P}(B_x(A_{n-1} \cup A_n)\vert C).$
Here we can already totally eliminate the $C$-contextuality for $B_x:$
\[{\bf P}(B_x(A_{n-1} \cup A_n)\vert C)={\bf P}(B_x\vert A_{n-1}) {\bf P} (A_{n-1}\vert C)+
{\bf P}(B_x\vert A_{n}) {\bf P} (A_{n}\vert C) \]
\[+2\lambda (B_x\vert \{ A_{n-1}, A_n\})
\sqrt{{\bf P}(B_x\vert A_{n-1}) {\bf P}(A_{n-1}\vert C) {\bf P}(B_x\vert A_n) {\bf P}(A_n\vert C)},\]
where the coefficient of statistical disturbance $\lambda$ was defined by (\ref{LA}).
And if $|\lambda|$ is bounded by 1 then we can represent the probability as the square
of the absolute value of the complex amplitude:
$$
\varphi_C^{(n-1)}(x)=\sqrt{{\bf P}(B_x\vert A_{n-1}) {\bf P}(A_{n-1}\vert C)}  +
e^{i \theta_C(x)} \sqrt{{\bf P}(B_x\vert A_n) {\bf P}(A_n\vert C)},
$$
where $\theta_C(x)=\pm \arccos \lambda (x\vert \{A_{n-1}, A_n\}, C).$

We have:
$$
\varphi_C^{(j)}(x)=\sqrt{{\bf P}(B_x(A_j \cup \ldots \cup A_n)\vert C)} \; e^{i \alpha_C^{(j)}(x)},
$$
where $\alpha_C^{(j)}(x)=\arg \varphi_C^{(j)}(x)
= \arccos \frac{M_j}{N_j},$
where $ M_j=\sqrt{{\bf P}(B_x\vert A_j) {\bf P}(A_j\vert C)}$\\
$+\mu (B_x\vert \{A_j, A_{j+1} \cup \ldots \cup A_n\}, C)\sqrt{{\bf P}(B_x(A_{j+1}\cup \ldots \cup A_n)\vert C)},$\\
$
N_j=\sqrt{{\bf P}(B_x (A_j \cup \ldots \cup A_n) \vert C)}.
$
Finally, we have:
$$
\alpha_C^{(n-1)}(x)=\arg \varphi_C^{(n-1)}(x)
$$
$$
=\arccos \frac{\sqrt{{\bf P}(B_x\vert A_{n-1}) {\bf P}(A_{n-1}\vert C)} +
\lambda (B_x\vert \{A_{n-1}, A_n\},C) \sqrt{{\bf P}(B_x\vert A_n) {\bf P}(A_n\vert C)}}
{\sqrt{{\bf P}(B_x(A_{n-1} \cup A_n)\vert C)}}.
$$
Thus we have:
$$
\varphi_C(x)=\sqrt{{\bf P}(B_x\vert A_1) {\bf P}(A_1\vert C)} + e^{i [\gamma_C^{(1)}(x) - \alpha_C^{(2)}(x)]}
\varphi_C^{(2)}(x)
$$
$$
=\sqrt{{\bf P}(B_x\vert A_1) {\bf P}(A_1\vert C)} + e^{i \beta_C^{(2)}(x)} \sqrt{{\bf P}(B_x\vert A_2) {\bf P}(A_2\vert C)}
$$
$$
+ e^{i\beta_C^{(3)}(x)} \varphi_C^{(3)}(x),
$$
where
$$
\beta_C^{(2)}(x)=\gamma_C^{(1)}(x) - \alpha_C^{(2)}(x), \beta_C^{(3)}(x)=\beta_C^2(x) +
\gamma_C^{(2)}(x) - \alpha_C^{(3)}(x).
$$
Finally, we obtain:
$$
\varphi_C(x)=\sum_{j=1}^n e^{i\beta_C^{(j)}(x)} \sqrt{{\bf P}(B_x\vert A_j) {\bf P}(A_j\vert C)}
$$
with $\beta_C^{(1)}(x)=0$ (this is just due to our special choice of a representation)
and $\beta_C^{(n)}(x)=\beta_C^{(n-1)}(x) + \theta_C (x).$

Thus by inductive splitting of multivalued variables into dichotomous
variables we represented contextual probabilities by complex amplitudes
$\varphi_C(x).$ Here the Born's rule holds true.

By using the standard in this paper
symbols
$p(x\vert y)={\bf P}(B_x\vert A_y)$ and $p_C^b(x)={\bf P}(B_x\vert C), p_C^a(y)={\bf P}(A_y\vert C)$
we write
$$
\varphi_C(x)=\sum_y e^{i \beta_C^{(y)}(x)} \sqrt{p_C^a(y) p(x\vert y)}.
$$
In particular, for $n=3$ we have
$$
\varphi_C(x)=\sqrt{p_C^a(a_1) p(x\vert a_1)} + e^{i \beta_C^{(2)} (x)} \sqrt{p_C^a(a_2) p(x\vert a_2)},
+ e^{i \beta_C^{(3)}(x)} \sqrt{p_C^a(a_3) p(x\vert a_3)} ,
$$
where
$$
\beta_C^{(2)}(x)=\gamma_C^{(1)}(x) - \alpha_C^{(2)}(x),
\beta_C^{(3)}(x)=\beta_C^{(2)}(x)+ \theta_C(x).
$$
We remark that each phase
$\beta_C^{(j)}(x)$ depends on all three $a$-contexts, $A_1, A_2, A_3.$ So
we cannot use the symbol $\beta_C(x\vert y).$ In $\beta_C^{(y)}(x)$ the $y$ is just the summation index;
in fact, $\beta_C^{(y)}(x)\equiv \beta_C^{(y)}(x\vert A_1, A_2, A_3).$ We remark that the probability
$p_C^b(x)$ can be represented as
$$
p_C^b(x)=|\varphi_C(x)|^2=\sum_y p_C^a(y) p(x\vert y)
$$
$$
+ 2\sum_{y_1<y_2}\cos [\beta_C^{(y_2)}(x)-\beta_C^{(y_1)}(x)] \sqrt{p_C^a (y_1) p_C^a (y_2) p(x\vert y_1) p(x\vert y_2)} .
$$

We can proceed in the same way as in the case of dichotomous random variables, see section 3,4.

I would like to thank A. V. Bulinskii, A. N. Shirayev, A. S. Holevo and V.M. Maximov
for discussions on probabilistic foundations.

{\bf References}

1. A. N. Kolmogoroff,  {\it Grundbegriffe der Wahrscheinlichkeitsrech,}
 Springer Verlag,  Berlin, 1933; reprinted :
{\it Foundations of the Probability Theory,}  Chelsea Publ. Comp., New York, 1956.

2.  A. N. Shiryayev, {\it   Probability,}  Springer, New York-Berlin-Heidelberg, 1984.

3.  A. S. Holevo, {\it Statistical structure of quantum theory,}  Springer, New York-Berlin-Heidelberg,
  2001.

4.  A. Yu. Khrennikov,  {\it Non-Kolmogorovian theories of probability and quantum physics,}
 Nauka, Fizmatlit, Moscow 2003 (in Russian).

5.   A. V.  Bulinskii and A. Yu. Khrennikov,   Nonclassical total probability formula and quantum
interference probabilities,  Statistics and Probability Letters, 70  (2004), pp. 49-58.

6.  R. Feynman and A. Hibbs, {\it Quantum Mechanics and Path
Integrals,} McGraw-Hill, New-York, 1965.

7. P. A. M.  Dirac, {\it The Principles of Quantum Mechanics,}
Oxford Univ. Press, Oxford, 1930.

8. J. von Neumann, {\it Mathematical foundations of quantum
mechanics,} Princeton Univ. Press, Princeton, N.J., 1955.

9.   A. N. Kolmogorov,  {\em The Theory of Probability.} In:  A. D. Alexandrov,
 A. N. Kolmogorov,  M. A. Lavrent'ev, eds. {\em Mathematics, Its Content, Methods,
and Meaning} 2, M.I.T. Press, Boston,  1965, pp. 110-118.

10. B. V.  Gnedenko, {\em The theory of probability,} Chelsea Publ. Com.,
New-York, 1962.

11.  A. Renye,  On a new axiomatics of probability theory, {\em Acta Mat. Acad. Sc. Hung.},
6 (1955), pp. 285-335.

12.  R. Von Mises, {\em The mathematical theory of probability and
 statistics,} Academic, London, 1964.

13.  E. Conte, O. Todarello,  A. Federici, F.  Vitiello, M. Lopane, A. Yu.  Khrennikov,
A preliminar evidence of quantum-like behaviour in measurements of mental states,
  Proc. Int. Conf. {\em Quantum Theory: Reconsideration
of Foundations,} A. Yu. Khrennikov , ed, Ser. Math. Modelling, 10,
 V\"axj\"o Univ. Press, V\"axj\"o, 2004, pp. 441-454.

14. A. Yu.  Khrennikov,  Quantum-like formalism for cognitive measurements, {\em Biosystems,}
70 (2003), pp. 211-233.

15. A. Yu. Khrennikov,{\it Information dynamics in cognitive, psychological and anomalous phenomena,}
Kluwer, Dordreht, 2004.

16. A. A. Grib, A. Yu. Khrennikov, K. Starkov, Probability amplitude in quantum-like games,
Proc. Int. Conf. {\it Quantum Theory: Reconsideration
of Foundations,} A. Yu. Khrennikov , ed., Ser. Math. Modelling, 10,
 V\"axj\"o Univ. Press,  V\"axj\"o , 2004, pp. 703-722.

17. O. A. Choustova, Bohmian mechanics for financial processes,
{\it J. Modern Optics}, 51 (2004), pp. 1111-1112.

18. A. Yu.  Khrennikov and S. V. Kozyrev,  Noncommutative probability in classical disordered systems,
{\em Physica A,} 326 (2003), pp. 456-463.

19. A. Yu.  Khrennikov,  Linear representations of probabilistic transformations
induced by context transitions, {\em J. Phys.A: Math. Gen.,} 34 (2001), pp. 9965-9981.

20. A. Yu.  Khrennikov,   Contextual viewpoint to quantum statistics,
{\em J. Math. Phys.}, 44 (2003), pp. 2471- 2478.

21.  A. Yu. Khrennikov, Interference of probabilities and number field structure of quantum models,
{\it Annalen  der Physik,} 12 (2003), pp. 575-585.

22. V. M. Maximov, Abstract models of Probability,  Proc. Conf.
{\it Foundations of Probability and Physics,} {\it Quantum Probability and White Noise Analysis},
{\bf 13}, 201- 218, WSP, Singapore, 2001, pp.257-273.
\end{document}